\documentclass[twoside,10pt]{article}
\usepackage{amssymb,amsmath,amsfonts,graphicx,rawfonts,t1enc}
\input{prepictex}
\input{pictex}
\input{postpictex}
\newcommand{\wt}{\widetilde}
\newcommand{\R}{\mathbb R}
\newcommand{\C}{\mathbb C}
\newcommand{\Q}{\mathbb Q}
\newcommand{\F}{\mathbb F}
\newcommand{\K}{\mathbb K}
\begin{document}
\thispagestyle{plain}
\hfill
\vskip 1.5truecm
\par
\noindent
{\large \bf A discrete form of the Beckman-Quarles theorem}
\par
\noindent
{\large \bf for mappings from ${\R}^2$ (${\C}^2$) to ${\F}^2$, where $\F$ is}
\par
\noindent
{\large \bf a subfield of a commutative field extending $\R$ ($\C$)}
\vskip 0.5truecm
\par
\noindent
{\large Apoloniusz Tyszka}
\vskip 0.5truecm
\par
\noindent
{\it Abstract.} Let $\F$ be a subfield of a commutative
field extending $\R$. Let
$\varphi_n:{\F}^n \times {\F}^n \to \F$,
$\varphi_n((x_1,...,x_n),(y_1,...,y_n))=(x_1-y_1)^2+...+(x_n-y_n)^2$.
We say that $f:{\R}^n \to {\F}^n$ preserves distance $d \geq 0$ if
$|x-y|=d$ implies $\varphi_n(f(x),f(y))=d^2$
for each $x,y \in {\R}^n$.
Let $A_{n}(\F)$ denote the set of all
positive numbers $d$ such that any map $f:{\R}^n \to {\F}^n$
that preserves unit distance preserves also distance $d$.
Let $D_{n}(\F)$ denote the set of all positive numbers~$d$
with the property:
if $x,y \in {\R}^n$ and $|x-y|=d$ then there exists a finite set $S_{xy}$
with $\left\{x,y \right\} \subseteq S_{xy} \subseteq {\R}^n$ such that any map
$f:S_{xy} \to {\F}^n$ that preserves unit distance
satisfies $|x-y|^2=\varphi_n(f(y),f(y))$.
Obviously, $\{1\} \subseteq D_n(\F) \subseteq A_{n}(\F)$.
We prove:
$A_{n}(\C) \subseteq \left\{d>0: d^2 \in \Q \right\} \subseteq D_2(\F)$.
Let $\K$ be a subfield of a commutative field $\Gamma$ extending $\C$.
Let
$\psi_2:{\Gamma}^2 \times {\Gamma}^2 \to \Gamma$,
$\psi_2((x_1,x_2),(y_1,y_2))=(x_1-y_1)^2+(x_2-y_2)^2$.
We say that $f:{\C}^2 \to {\K}^2$ preserves unit distance if
$\psi_2(X,Y)=1$ implies $\psi_2(f(X),f(Y))=1$
for each $X,Y \in {\C}^2$.
We prove: if $X,Y \in {\C}^2$, $\psi_2(X,Y) \in \Q$
and $X \neq Y$, then there exists a finite set $S_{XY}$ with
$\left\{X,Y \right\} \subseteq S_{XY} \subseteq {\C}^2$ such that
any map $f:S_{XY} \to {\K}^2$ that preserves unit distance
satisfies $\psi_2(X,Y)=\psi_2(f(X),f(Y))$ and $f(X) \neq f(Y)$.
\vskip 0.5truecm
\par
\noindent
Let $\F$ be a subfield of commutative field extending $\R$.
Let $\varphi_n: {\F}^n \times {\F}^n \to \F$,
$\varphi_n((x_1,...,x_n),(y_1,...,y_n))=(x_1-y_1)^2+...+(x_n-y_n)^2$.
We say that $f:{\R}^n \to {\F}^n$ preserves distance $d \geq 0$ if
$|x-y|=d$ implies $\varphi_n(f(x),f(y))=d^2$ for each $x,y \in {\R}^n$.
We say that $f:{\R}^n \to {\F}^n$ preserves the distance between
$x,y \in {\R}^n$ if $|x-y|^2=\varphi_n(f(x),f(y))$.
By a {\sl complex isometry} of ${\C}^n$ (i.e. an affine isometry of ${\C}^n$)
we understand any map $h:{\C}^n \to {\C}^n$ of the form
$$
h(z_1,z_2,...,z_n)=(z_1',z_2',...,z_n')
$$
where
$$
z_j'=a_{0j}+a_{1j}z_{1}+a_{2j}z_{2}+...+a_{nj}z_{n} \hspace{0.5truecm}(j=1,2,...,n),
$$
the coefficients $a_{ij}$ are complex and the matrix $||a_{ij}||$
$(i,j=1,2,...,n)$ is orthogonal i.e. satisfies the condition
$$
\sum_{j=1}^{n} a_{\mu j} a_{\nu j}=\delta_{\nu}^{\mu} \hspace{0.5truecm} (\mu,\nu=1,2,...,n)
$$
with Kronecker's delta. According to \cite{Borsuk},
$\varphi_n(x,y)$ is invariant under complex isometries i.e. for every
complex isometry $h:{\C}^n \to {\C}^n$
\begin{equation}
\tag*{{\boldmath $\left(\diamond\right)$}}
\forall x,y \in {\C}^n ~~\varphi_n(h(x),h(y))=\varphi_n(x,y)
\end{equation}
Conversely, if $h:{\C}^n \to {\C}^n$ satisfies
{\boldmath $\left(\diamond\right)$} then $h$ is a complex isometry;
it follows from [4, proposition 1, page 21] by replacing $\R$ with $\C$
and $d(x,y)$ with $\varphi_n(x,y)$. Similarly, if
$f:{\R}^n \to {\C}^n$ preserves all distances, then there exists
a complex isometry $h:{\C}^n \to {\C}^n$ such that
$f=h_{|{\R}^n}$.
\vskip 0.1truecm
\par
\noindent
By a field endomorphism of a field $H$ we understand any ring
endomorphism $g:H \to H$ with $g(1) \neq 0$.
Bijective endomorphisms are called automorphisms,
for more information on field endomorphisms and automorphisms
of $\C$ the reader is referred to \cite{Kestelman}, \cite{Kuczma} and \cite{Yale}.
If $r$ is a rational number, then $g(r)=r$ for any field
endomorphism $g:\C \to \C$.
Proposition~1 shows that only rational numbers $r$ have this property:
\vskip 0.1truecm
\par
\noindent
{\bf Proposition 1} (\cite{Tyszka2004a}).
If $r \in \C$ and $r$ is not a rational number,
then there exists a field automorphism $g:\C \to \C$ such that
$g(r) \neq r$.
\vskip 0.1truecm
\par
\noindent
{\bf Theorem~1} (\cite{Tyszka2004a}). If $x,y \in {\R}^n$ and $|x-y|^2$ is not a rational
number, then there exists $f:{\R}^n \to {\C}^n$ that does not preserve the
distance between $x$ and~$y$ although $f$ preserves all distances
$\sqrt{r}$ with rational $r \geq 0$.
\vskip 0.1truecm
\par
\noindent
Let $A_{n}(\F)$ denote the set of all positive numbers $d$
such that any map $f:{\R}^n \to {\F}^n$ that preserves
unit distance preserves also distance $d$.
The classical Beckman-Quarles theorem states that each unit-distance
preserving mapping from ${\R}^n$ to ${\R}^n$ ($n \geq 2$) is an isometry,
see \cite{Beckman-Quarles}-\cite{Benz1994} and \cite{Everling}.
It means that $A_{n}(\R)=(0,\infty)$ for each $n \geq 2$.
By Theorem~1 $A_{n}(\C) \subseteq \left\{d>0: d^2 \in \Q \right\}$.
Let $D_{n}(\F)$ denote the set of all positive numbers $d$
with the following property:
\begin{description}
\item{{\boldmath $\left(\star\right)$}}
if $x,y \in {\R}^n$ and $|x-y|=d$ then there exists a finite set $S_{xy}$
with $\left\{x,y \right\} \subseteq S_{xy} \subseteq {\R}^n$ such that any map
$f:S_{xy}\to {\F}^n$ that preserves unit distance
preserves also the distance between $x$ and $y$.
\end{description}
\par
\noindent
Obviously, $\{1\} \subseteq D_n(\F) \subseteq A_{n}(\F)$.
The author proved in \cite{Tyszka2004a}
that $D_n(\C)$ is a dense subset of $(0,\infty)$ for each $n \geq 2$.
The proof remains valid for $D_n(\F)$.
It is also known~(\cite{Tyszka2000},~\cite{Tyszka2001}) that
$D_n(\R)$ is equal to the set of positive algebraic numbers
for each $n \geq 2$.
We shall prove that $\{d>0: d^2 \in \Q\} \subseteq D_2(\F)$.
We need the following technical Propositions 2-5.
\vskip 0.1truecm
\par
\noindent
{\bf Proposition~2} (cf. \cite{Blumenthal}, \cite{Borsuk}).
The points
$$c_{1}=(z_{1,1},~...,~z_{1,n}),~...,~c_{n+1}=(z_{n+1,1},~...,~z_{n+1,n}) \in {\F}^n$$
are affinely dependent if and only if their Cayley-Menger determinant
$$\Delta(c_1,...,c_{n+1}):=
\det \left[
\begin{array}{ccccc}
 0  &  1                       &  1                       & ... & 1                       \\
 1  &  0                       & \varphi_n(c_{1},c_{2})   & ... & \varphi_n(c_{1},c_{n+1})\\
 1  & \varphi_n(c_{2},c_{1})   &  0                       & ... & \varphi_n(c_{2},c_{n+1})\\
... & ...                      & ...  	                  & ... & ...                     \\
 1  & \varphi_n(c_{n+1},c_{1}) & \varphi_n(c_{n+1},c_{2}) & ... & 0                       \\
\end{array}\;\right]
$$
\par
\noindent
equals $0$.
\vskip 0.1truecm
\par
\noindent
{\it Proof.} It follows from the equality
$$
\left(
\det \left[
\begin{array}{ccccc}
z_{1,1}   & z_{1,2}   & ... &  z_{1,n}  & 1  \\
z_{2,1}   & z_{2,2}   & ... &  z_{2,n}  & 1  \\
  ...     &  ...      & ... &  ...      & ...\\
z_{n+1,1} & z_{n+1,2} & ... & z_{n+1,n} & 1  \\
\end{array}
\right] \right)^2=
\frac{(-1)^{n+1}}{2^{n}} \cdot \Delta(c_1,...,c_{n+1}).
$$
\vskip 0.1truecm
\par
\noindent
{\bf Proposition~3} (cf. \cite{Blumenthal}, \cite{Borsuk}).
For any points $c_{1},..., c_{n+k} \in {\F}^n$ ($k=2,3,4,...$) their
Cayley-Menger determinant equals $0$ i.e. $\Delta(c_1,...,c_{n+k})=0$.
\vskip 0.1truecm
\par
\noindent
{\it Proof.} Assume that
$c_{1}=(z_{1,1},~...,~z_{1,n})$,~...,~~$c_{n+k}=(z_{n+k,1},~...,~z_{n+k,n})$.
The points
$$\wt{c}_{1}=(z_{1,1},~...,~z_{1,n},~0,...,0),~...,~\wt{c}_{n+k}=(z_{n+k,1},~...,~z_{n+k,n},~0,...,0) \in {\F}^{n+k-1}$$
are affinely dependent.
Since $\varphi_n(c_i,c_j)=\varphi_{n+k-1}(\wt{c}_{i},\wt{c}_{j})~(1 \leq i \leq j \leq n+k)$
the Cayley-Menger determinant
of points $c_{1},...,c_{n+k}$ is equal to the Cayley-Menger determinant
of points $\wt{c}_{1},...,\wt{c}_{n+k}$ which equals $0$ according to
Proposition~2.
\vskip 0.1truecm
\par
\noindent
{\bf Proposition~4a.}
If $d \in \F$, $d \neq 0$, $c_{1},c_{2},c_{3} \in {\F}^2$ and
$\varphi_2(c_1,c_2)=\varphi_2(c_1,c_3)=\varphi_2(c_2,c_3)=d^2$,
then $c_{1}$, $c_{2}$, $c_{3}$ are affinely independent.
\vskip 0.1truecm
\par
\noindent
{\it Proof.} It follows from Proposition 2 because
the Cayley-Menger determinant
$$\Delta(c_1,c_2,c_3)=
\det \left[
\begin{array}{cccc}
0 &  1  &  1  &  1  \\
1 &  0  & d^2 & d^2 \\
1 & d^2 &  0  & d^2 \\
1 & d^2 & d^2 &  0  \\
\end{array}
\right]=-3d^4 \neq 0.
$$
\par
\noindent
{\bf Proposition~4b} (\cite{Tyszka2004b}). If $a,b \in \F$, $a+b \neq 0$, $z,x,\wt{x} \in {\F}^2$
and $\varphi_2(z,x)=a^2$, $\varphi_2(x,\wt{x})=b^2$,
$\varphi_2(z,\wt{x})=(a+b)^2$, then
$\overrightarrow{zx}=\frac{\textstyle a}{{\textstyle a+b}}\overrightarrow{z\wt{x}}$.
\vskip 0.1truecm
\par
\noindent
{\bf Proposition~5} (cf. [6, Lemma, page 127] for ${\R}^n$).
If $x$, $y$, $c_{0}$, $c_{1}$, $c_{2}$ $\in$ ${\F}^2$,
$\varphi_2(x,c_{0})=\varphi_2(y,c_{0})$,
$\varphi_2(x,c_{1})=\varphi_2(y,c_{1})$,
$\varphi_2(x,c_{2})=\varphi_2(y,c_{2})$
and $c_0$, $c_1$, $c_2$ are affinely independent,
then $x=y$.
\vskip 0.1truecm
\par
\noindent
{\it Proof.} Computing we obtain that the vector
$\overrightarrow{xy}$ is perpendicular
to each of the linearly independent vectors
$\overrightarrow{c_{0}c_{1}}$, $\overrightarrow{c_{0}c_{2}}$.
Thus the vector $\overrightarrow{xy}$ is perpendicular
to every linear combination of vectors
$\overrightarrow{c_{0}c_{1}}$ and $\overrightarrow{c_{0}c_{2}}$.
In particular, the vector
$\overrightarrow{xy}$ is perpendicular
to each of the vectors $[1,0]$ and $[0,1]$.
Therefore $\overrightarrow{xy}=0$ and the proof is completed.
\vskip 0.1truecm
\par
\noindent
{\bf Lemma~1}. If $d \in D_2(\F)$ then $\sqrt{3} \cdot d \in D_2(\F)$.
\vskip 0.1truecm
\par
\noindent
{\it Proof.}
Let $ d \in D_2(\F)$, $x,y \in {\R}^2$ and $|x-y|=\sqrt{3} \cdot d$.
Using the notation of Figure~1 we show that
$$
S_{xy}:=S_{y\wt{y}} \cup \bigcup_{i=1}^2 S_{xp_i} \cup \bigcup_{i=1}^2 S_{yp_i} \cup S_{p_1p_2}
\cup
\bigcup_{i=1}^2 S_{x\wt{p_i}}
\cup
\bigcup_{i=1}^2 S_{\wt{y}\wt{p_i}}
\cup S_{\wt{p_1}\wt{p_2}}
$$
satisfies condition {\boldmath $\left(\star\right)$}.
\vskip 0.1truecm
\par
\centerline{
\beginpicture
\setcoordinatesystem units <0.5cm, 0.5cm>
\setplotsymbol({ .})
\plot 6.00 0.00 0.00 0.00 /
\plot 0.00 0.00 3.00 5.20 /
\plot 3.00 5.20 6.00 0.00 /
\plot 6.00 0.00 9.00 5.20 /
\plot 9.00 5.20 3.00 5.20 /
\plot 9.38 -0.80 4.38 -4.10 /
\plot 4.38 -4.10 4.00 1.88 /
\plot 9.0 5.20 9.38 -0.80 /
\plot 9.00 5.20 4.00 1.88 /
\plot 0.00 0.00 4.38 -4.10 /
\plot 9.38 -0.80 4.00 1.88 /
\setdots
\plot 0.00 0.00 9.0 5.2 /
\plot 9.0 5.20 4.38 -4.10 /
\put {$\wt{p_1}$} at 10 -0.80
\put {$y$} at -0.20 0.00
\put {$p_1$} at 5.80 -0.40
\put {$x$} at 9.00 5.60
\put {$p_2$} at 3.00 5.60
\put {$\wt{y}$} at 4.38 -4.60
\put {$\wt{p_2}$} at 3.60 1.40
\endpicture}
\vskip 0.1truecm
\par
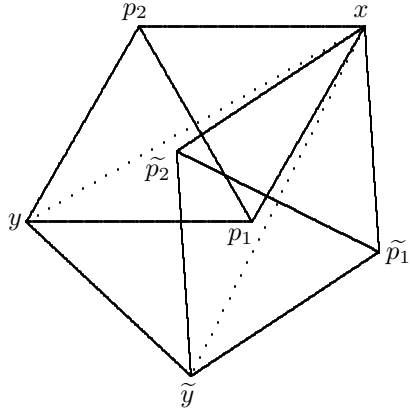
\centerline{Figure~1}
\centerline{$|x-y|=|x-\wt{y}|=\sqrt{3} \cdot d$, $|y-\wt{y}|=|p_1-p_2|=|\wt{p_1}-\wt{p_2}|=d$}
\centerline{$|x-p_i|=|y-p_i|=|x-\wt{p_i}|=|\wt{y}-\wt{p_i}|=d$ ($i=1,2$)}
\vskip 0.1truecm
\par
\noindent
Assume that $f: S_{xy} \to {\F}^2$
preserves unit distance. Since
$$
S_{xy} \supseteq
S_{y\wt{y}}
\cup
\bigcup_{i=1}^{2}S_{xp_{i}}
\cup\bigcup_ {i=1}^{2}S_{yp_{i}}
\cup
S_{p_{1}p_{2}}
$$
\noindent
we conclude that $f$ preserves the distances between
$y$ and $\wt{y}$,
$x$ and $p_{i}$ ($i=1,2$), $y$ and $p_{i}$ ($i=1,2$),
$p_{1}$ and $p_{2}$.
Hence $\varphi_2(f(y),f(\wt{y}))=
\varphi_2(f(x),f(p_i))=\varphi_2(f(y),f(p_i))=\varphi_2(f(p_1),f(p_2))=d^2$
($i=1,2$).
By Proposition~3 the Cayley-Menger determinant
$\Delta(f(x),f(p_{1}),f(p_{2}),f(y))$ equals $0$ i.e.
$$
\det\left[
\begin{array}{ccccc}
0 &  1                     &  1                       &  1                       &  1                    \\
1 &  0                     & \varphi_2(f(x),f(p_1))   & \varphi_2(f(x),f(p_2))   & \varphi_2(f(x),f(y))  \\
1 & \varphi_2(f(p_1),f(x)) &  0                       & \varphi_2(f(p_1),f(p_2)) & \varphi_2(f(p_1),f(y))\\
1 & \varphi_2(f(p_2),f(x)) & \varphi_2(f(p_2),f(p_1)) &  0                       & \varphi_2(f(p_2),f(y))\\
1 & \varphi_2(f(y),f(x))   & \varphi_2(f(y),f(p_1))   & \varphi_2(f(y),f(p_2))   &  0                    \\
\end{array}
\right]
=0.
$$
Therefore
$$
\det \left[
\begin{array}{ccccccc}
0 &  1  &  1  &  1  &  1 \\
1 &  0  & d^2 & d^2 &  t \\
1 & d^2 &  0  & d^2 & d^2\\
1 & d^2 & d^2 &  0  & d^2\\
1 &  t  & d^2 & d^2 &  0 \\
\end{array}
\right]
=0
$$
where $t=\varphi_2(f(x),f(y))$. Computing this determinant we obtain
\par
\noindent
\centerline{$2d^{2}t \cdot (3d^2-t)=0.$}
\par
\noindent
Therefore
\par
\noindent
\centerline{$t=\varphi_2(f(x),f(y))=\varphi_2(f(y),f(x))=3d^2=|x-y|^2$}
or
\par
\noindent
\centerline{$t=\varphi_2(f(x),f(y))=\varphi_2(f(y),f(x))=0.$}
\par
\noindent
Analogously we may prove that
\vskip 0.1truecm
\par
\noindent
\centerline{$\varphi_2(f(x),f(\wt{y}))=\varphi_2(f(\wt{y}),f(x))=3d^2=|x-\wt{y}|^2$}
\par
\noindent
or
\par
\noindent
\centerline{$\varphi_2(f(x),f(\wt{y}))=\varphi_2(f(\wt{y}),f(x))=0.$}
\vskip 0.1truecm
\par
\noindent
If $t=0$ then the points $f(x)$ and $f(y)$ satisfy:
\par
\noindent
\centerline{$\varphi_2(f(x),f(x))=0=\varphi_2(f(y),f(x))$,}
\vskip 0.1truecm
\centerline{$\varphi_2(f(x),f(p_1))=d^2=\varphi_2(f(y),f(p_1))$,}
\vskip 0.1truecm
\centerline{$\varphi_2(f(x),f(p_2))=d^2=\varphi_2(f(y),f(p_2))$.}
\vskip 0.1truecm
\par
\noindent
By Proposition~4a the points $f(x)$, $f(p_1)$, $f(p_2)$ are affinely
independent. Therefore by Proposition~5 $f(x)=f(y)$ and consequently
\par
\noindent
\centerline{$d^2=\varphi_2(f(y),f(\wt{y}))=\varphi_2(f(x),f(\wt{y}))
\in \{3d^2,~0 \}$.}
\par
\noindent
Since $d^2 \neq 3d^2$ and $d^2 \neq 0$ we conclude that
the case $t=0$ cannot occur. This completes the proof of Lemma~1.
\vskip 0.1truecm
\par
\noindent
{\bf Lemma~2.} If $d \in D_2(\F)$ then $2 \cdot d \in D_2(\F)$.
\vskip 0.1truecm
\par
\noindent
{\it Proof.}
Let $d \in D_2(\F)$, $x,y \in {\R}^2$ and $|x-y|=2 \cdot d$.
Using the notation of Figure~2 we show that
\vskip 0.1truecm
\centerline{$S_{xy}:=
\bigcup \{S_{ab}:a,b\in \{x,y,p_{1},p_{2},p_{3}\},
|a-b|=d \vee |a-b|=\sqrt{3} \cdot d\}$}
\vskip 0.1truecm
\par
\noindent
(where $S_{xp_3}$ and $S_{yp_2}$ are known to exist by Lemma~1)
satisfies condition {\boldmath $\left(\star\right)$}.
\par
\centerline{
\beginpicture
\normalsize
\setcoordinatesystem units <1mm,1mm>
\setplotarea x from -30 to 30, y from -1 to 27
\setplotsymbol({ .})
\plot 0 0 30 0 /
\plot 0 0 -30 0 /
\plot 30 0 15 26 /
\plot 15 26   0 0 /
\plot 0 0 -15 26 /
\plot -15 26 -30 0 /
\plot -15 26 15 26 /
\setdashes
\plot -30 0 15 26 /
\plot 30 0 -15 26 /
\put{$x$} at -30 -2
\put {$p_1$} at 1 -2
\put {$p_2$} at -15 28
\put{$p_3$} at 15 28
\put {$y$} at 30 -2
\endpicture}
\vskip 0.1truecm
\par
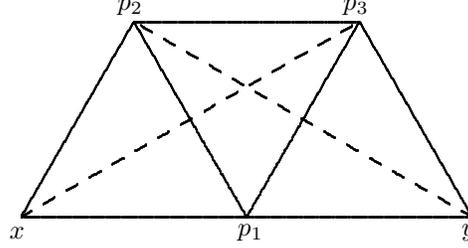
\centerline{Figure~2}
\centerline{$|x-y|=2 \cdot d$}
\centerline{$|p_1-p_2|=|p_1-p_3|=|p_2-p_3|=|x-p_1|=|x-p_2|=|y-p_1|=|y-p_3|=d$}
\centerline{$|x-p_3|=|y-p_2|=\sqrt{3} \cdot d$}
\vskip 0.1truecm
\par
\noindent
Assume that
$f:S_{xy} \to {\F}^2$ preserves unit distance. Then
$f$ preserves all distances between $p_i$ and $p_j$ $(1 \leq i < j \leq 3)$,
$x$ and $p_i$ $(1 \leq i \leq 3)$, $y$ and $p_i$ $(1 \leq i \leq 3)$.
By Proposition~3 the Cayley-Menger determinant
$\Delta(f(x)$, $f(p_{1})$, $f(p_{2})$, $f(p_{3})$, $f(y))$ equals $0$ i.e.
\arraycolsep=0.5mm
\footnotesize
$$
\det \left[
\begin{array}{cccccc}
0 &  1                     &  1                       &  1                       &  1                       &  1                    \\
1 &  0                     & \varphi_2(f(x),f(p_1))   & \varphi_2(f(x),f(p_2))   & \varphi_2(f(x),f(p_3))   & \varphi_2(f(x),f(y))  \\
1 & \varphi_2(f(p_1),f(x)) &  0                       & \varphi_2(f(p_1),f(p_2)) & \varphi_2(f(p_1),f(p_3)) & \varphi_2(f(p_1),f(y))\\
1 & \varphi_2(f(p_2),f(x)) & \varphi_2(f(p_2),f(p_1)) &  0                       & \varphi_2(f(p_2),f(p_3)) & \varphi_2(f(p_2),f(y))\\
1 & \varphi_2(f(p_3),f(x)) & \varphi_2(f(p_3),f(p_1)) & \varphi_2(f(p_3),f(p_2)) &  0                       & \varphi_2(f(p_3),f(y))\\
1 & \varphi_2(f(y),f(x))   & \varphi_2(f(y),f(p_1))   & \varphi_2(f(y),f(p_2))   & \varphi_2(f(y),f(p_3))   &  0                    \\
\end{array}
\right]
=0.
$$
\normalsize
\arraycolsep=1.0mm
Therefore
$$
\det \left[
\begin{array}{cccccc}
0 &  1   &  1  &  1   &  1   &  1  \\
1 &  0   & d^2 & d^2  & 3d^2 &  t  \\
1 & d^2  &  0  & d^2  & d^2  & d^2 \\
1 & d^2  & d^2 &  0   & d^2  & 3d^2\\
1 & 3d^2 & d^2 & d^2  &  0   & d^2 \\
1 &  t   & d^2 & 3d^2 & d^2  &  0  \\
\end{array}
\right]
=0
$$
where $t=\varphi_2(f(x),f(y))$. Computing this determinant we obtain
$$3d^4 \cdot (t-4d^2)^{2}=0.$$
Therefore
$$t=\varphi_2(f(x),f(y))=\varphi_2(f(y),f(x))=4d^2=|x-y|^2.$$
\par
\noindent
{\bf Lemma~3.} If $a,b \in D_2(\F)$ and $a>b$, then $\sqrt{a^2-b^2} \in D_2(\F)$.
\vskip 0.1truecm
\par
\noindent
{\it Proof.}
Let $a, b \in D_2(\F)$, $a>b$, $x,y \in {\R}^2$ and $|x-y|=\sqrt{a^2-b^2}$. Using the notation
of Figure~3 we show that
$$
S_{xy}:=S_{xp_1} \cup S_{xp_2} \cup S_{yp_1} \cup S_{yp_2} \cup S_{p_1p_2}
$$
(where $S_{p_1p_2}$ is known to exist by Lemma~2)
satisfies condition {\boldmath $\left(\star\right)$}.
\vskip 0.1truecm
\par
\centerline{
\beginpicture
\setcoordinatesystem units <1mm, 1mm>
\setplotarea x from -27 to 27, y from -2 to 44
\setplotsymbol({ .})
\plot -25 0 0 0 /
\plot 0 0 25 0 /
\plot 25 0 0 42 /
\plot -25 0 0 42 /
\setdots
\plot 0 0 0 42 /
\put {$x$} at 0.5 -2
\put {$y$} at 0.5 44
\put {$p_1$} at -25 -2
\put{$p_2$} at 25 -2
\put {$b$} at -12.5 2.5
\put {$b$} at  12.5 2.5
\put{$a$} at -14 22
\put{$a$} at 15 23
\endpicture}
\vskip 0.1truecm
\par
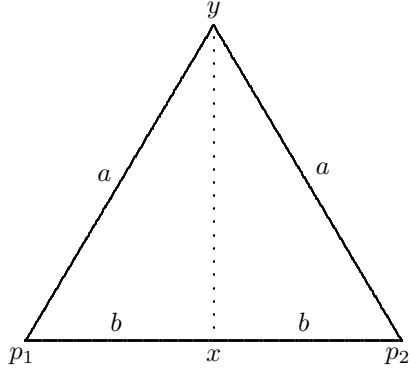
\centerline{Figure~3}
\centerline{$|x-y|=\sqrt{a^2-b^2}$}
\centerline{$|x-p_1|=|x-p_2|=b$, $|y-p_1|=|y-p_2|=a$, $|p_1-p_2|=2b$}
\vskip 0.1truecm
\par
\noindent
Assume that
$f: S_{xy} \to {\F}^2$ preserves unit distance.
Then $f$ preserves the distances between
$x$ and $p_{i}$ ($i=1,2$), $y$ and $p_{i}$ ($i=1,2$),
$p_1$ and $p_2$.
By Proposition~3 the Cayley-Menger determinant
$\Delta(f(x)$, $f(p_1)$, $f(p_2)$, $f(y))$ equals $0$ i.e.
$$
\det \left[
\begin{array}{ccccc}
0 &  1                     &  1                       &  1                       &  1                    \\
1 &  0                     & \varphi_2(f(x),f(p_1))   & \varphi_2(f(x),f(p_2))   & \varphi_2(f(x),f(y))  \\
1 & \varphi_2(f(p_1),f(x)) &  0                       & \varphi_2(f(p_1),f(p_2)) & \varphi_2(f(p_1),f(y))\\
1 & \varphi_2(f(p_2),f(x)) & \varphi_2(f(p_2),f(p_1)) &  0                       & \varphi_2(f(p_2),f(y))\\
1 & \varphi_2(f(y),f(x))   & \varphi_2(f(y),f(p_1))   & \varphi_2(f(y),f(p_2))   &  0                    \\
\end{array}
\right]
=0.
$$
Therefore
$$
\det \left[
\begin{array}{ccccc}
0 &  1  &  1   &  1   &  1 \\
1 &  0  & b^2  & b^2  &  t \\
1 & b^2 &  0   & 4b^2 & a^2\\
1 & b^2 & 4b^2 &  0   & a^2\\
1 &  t  & a^2  & a^2  &  0 \\
\end{array}
\right]
=0
$$
where $t=\varphi_2(f(x),f(y))$. Computing this determinant we obtain
$$
-8b^2 \cdot (t+b^2-a^2)^2=0.
$$
Therefore $$t=a^2-b^2=|x-y|^2.$$
\vskip 0.1truecm
\par
\noindent
{\bf Lemma~4.} All distances $\sqrt{n}$~($n=1,2,3,...$) belong
to $D_2(\F)$.
\vskip 0.1truecm
\par
\noindent
{\it Proof.} If $n \in \{2,3,4,...\}$ and $\sqrt{n} \in D_2(\F)$,
then $\sqrt{n-1}=\sqrt{(\sqrt{n})^2-1^2} \in D_2(\F)$;
it follows from Lemma 3 and $1 \in D_2(\F)$.
On the other hand, by Lemma 2 all numbers $\sqrt{2^{2k}}=2^k$ ($k=0,1,2,..$) belong
to $D_2(\F)$. These two facts imply that all distances $\sqrt{n}$ ($n=1,2,3,...$) belong
to $D_2(\F)$.
\vskip 0.1truecm
\par
\noindent
From Lemma 4 we obtain:
\vskip 0.1truecm
\par
\noindent
{\bf Lemma 5.} For each $n \in \{1,2,3,...\}$ we have: $n=\sqrt{n^2} \in D_2(\F)$.
\vskip 0.1truecm
\par
\noindent
{\bf Lemma 6.} If $d \in D_2(\F)$ then all distances $\frac{\textstyle d}{\textstyle k}$
($k=2,3,4,...$) belong to $D_2(\F)$.
\vskip 0.1truecm
\par
\noindent
{\it Proof.} Let $d \in D_2(\F)$, $k \in \{2,3,4,...\}$, $x,y \in {\R}^2$ and
$|x-y|=\frac{\textstyle d}{\textstyle k}$. We choose an
integer $e \geq d$. Using the notation of Figure 4 we show that
$$
S_{xy}:=S_{\wt{x}\wt{y}}
\cup
S_{zx}
\cup
S_{x\wt{x}}
\cup
S_{z\wt{x}}
\cup
S_{zy}
\cup
S_{y\wt{y}}
\cup
S_{z\wt{y}}
$$
\par
\noindent
(where sets $S_{zx}$, $S_{x\wt{x}}$, $S_{z\wt{x}}$,
$S_{zy}$, $S_{y\wt{y}}$, $S_{z\wt{y}}$ corresponding
to integer distances are known
to exist by Lemma 5) satisfies condition {\boldmath $\left(\star\right)$}.
\vskip 0.1truecm
\par
\centerline{
\beginpicture
\setcoordinatesystem units <1mm, 1mm>
\setplotsymbol({ .})
\plot 0 0 -59.3211 9 /
\plot 0 0 -59.3211 -9 /
\plot -59.3211 -9 -59.3211 9 /
\setdots
\plot -19.7737 -3 -19.7737 3 /
\put {$z$} at 2 0
\put {$\wt{y}$} at -59.3211 11.5
\put {$\wt{x}$} at -59.3211 -11.5
\put {$y$} at -19.7737 5
\put {$x$} at -19.7737 -5
\put{$d$} at -61 0
\put {(} at -45 10.70
\put {$k$} at -42 10.25
\put {$-$} at -39 9.80
\put {$1$} at -36 9.35
\put {$)$} at -33 8.90
\put {$e$} at -30 8.45
\put {$e$} at -10 5.45
\put {(} at -45 -10.70
\put {$k$} at -42 -10.25
\put {$-$} at -39 -9.80
\put {$1$} at -36 -9.35
\put {$)$} at -33 -8.90
\put {$e$} at -30 -8.45
\put {$e$} at -10 -5.45
\endpicture}
\vskip 0.1truecm
\par
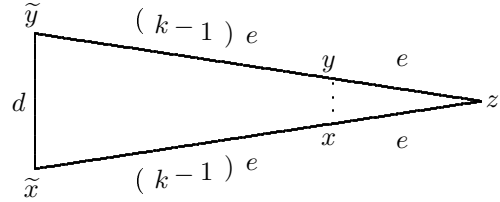
\centerline{Figure 4}
\centerline{$|x-y|=\frac{\textstyle d}{\textstyle k}$}
\vskip 0.1truecm
\par
\noindent
Assume that $f:S_{xy} \to {\F}^2$ preserves unit distance.
Since
$
S_{xy} \supseteq
S_{zx}
\cup
S_{x\wt{x}}
\cup
S_{z\wt{x}}
$
we conclude that:
\par
\noindent
\centerline{$\varphi_2(f(z),f(x))=|z-x|^2=e^2$,}
\\
\centerline{$\varphi_2(f(x),f(\wt{x}))=|x-\wt{x}|^2=((k-1)e)^2$,}
\\
\centerline{$\varphi_2(f(z),f(\wt{x}))=|z-\wt{x}|^2=(ke)^2$.}
\par
\noindent
By Proposition 4b:
\begin{equation}
\overrightarrow{f(z)f(x)}=\frac{1}{k}\overrightarrow{f(z)f(\wt{x})}
\end{equation}
Analogously:
\begin{equation}
\overrightarrow{f(z)f(y)}=\frac{1}{k}\overrightarrow{f(z)f(\wt{y})}
\end{equation}
By (1) and (2):
\begin{equation}
\overrightarrow{f(x)f(y)}=\frac{\textstyle 1}{\textstyle k}\overrightarrow{f(\wt{x})f(\wt{y})}
\end{equation}
Since
$S_{xy} \supseteq S_{\wt{x}\wt{y}}$
we conclude that
\begin{equation}
\varphi_2(f(\wt{x}),f(\wt{y}))=|\wt{x}-\wt{y}|^2=d^2
\end{equation}
By (3) and (4): 
$\varphi_2(f(x),f(y))=\left(\frac{\textstyle d}{\textstyle k}\right)^2$
and the proof is completed.
\vskip 0.1truecm
\par
\noindent
{\bf Theorem 2.} If $x,y \in {\R}^2$ and $|x-y|^2$ is a rational
number, then there exists a finite set $S_{xy}$ with
$\left\{x,y \right\} \subseteq S_{xy} \subseteq {\R}^2$ such that
any map $f:S_{xy}\to {\F}^2$ that preserves unit distance
preserves also the distance between $x$ and $y$; in other words
$\{d>0: d^2 \in \Q \} \subseteq D_2(\F)$.
\vskip 0.1truecm
\par
\noindent
{\it Proof.} We need to prove that
$\sqrt{\frac{\textstyle p}{\textstyle q}} \in D_2(\F)$ for all positive integers $p,q$.
Since $\sqrt{\frac{\textstyle p}{\textstyle q}}=\frac{\sqrt{\textstyle pq}}{\textstyle q}$
the assertion follows from Lemmas 4 and 6.
\vskip 0.3truecm
\par
\noindent
Since $D_2(\C) \subseteq A_2(\C) \subseteq \{d>0: d^2  \in \Q \}$
as a corollary we get:
\begin{equation}
\tag*{{\boldmath $\left(\bullet\right)$}}
D_2(\C)=A_2(\C)=\{d>0: d^2 \in \Q \}
\end{equation}
\vskip 0.1truecm
\par
Let $\K$ be a subfield of a commutative field $\Gamma$ extending $\C$.
Let $\psi_2:{\Gamma}^2 \times {\Gamma}^2 \to \Gamma$,
$\psi_2((x_1,x_2),(y_1,y_2))=(x_1-y_1)^2+(x_2-y_2)^2$.
We say that $f:{\C}^2 \to {\K}^2$ preserves unit distance if
$\psi_2(X,Y)=1$ implies $\psi_2(f(X),f(Y))=1$
for each $X,Y \in {\C}^2$. We say that $f:{\C}^2 \to {\K}^2$
preserves the distance between $X,Y \in {\C}^2$ if
$\psi_2(X,Y)=\psi_2(f(X),f(Y))$.
\vskip 0.1truecm
\par
\noindent
{\bf Theorem 3.} If $X,Y \in {\C}^2$, $\psi_2(X,Y) \in \Q$
and $X \neq Y$,
then there exists a finite set $S_{XY}$ with
$\left\{X,Y \right\} \subseteq S_{XY} \subseteq {\C}^2$ such that
any map $f:S_{XY} \to {\K}^2$ that preserves unit distance
satisfies $\psi_2(X,Y)=\psi_2(f(X),f(Y))$ and $f(X) \neq f(Y)$.
\vskip 0.1truecm
\par
\noindent
{\it Proof.} If $\psi_2(X,Y)=\psi_2(f(X),f(Y))$, then
$\psi_2(X,Y) \neq 0$ implies $f(X) \neq f(Y)$.
The main part of the proof divides into three cases.
\vskip 0.1truecm
\par
\noindent
Case 1: $\psi_2(X,Y)>0$.
\vskip 0.1truecm
\par
\noindent
There exists a complex isometry of ${\C}^2$ 
that sends $(0,0)$ to $X$ and $(\sqrt{\psi_2(X,Y)},0)$ to~$Y$,
so without loss of generality we may assume that
$X=(0,0)$ and $Y=(\sqrt{\psi_2(X,Y)},0)$.
By Theorem 2 there exists a finite set $S_{XY}$ with
$\{X,Y\} \subseteq S_{XY} \subseteq {\R}^2 \subseteq {\C}^2$
such that each unit-distance preserving
mapping $f:S_{XY} \to {\K}^2$ satisfies
$\psi_2(X,Y)=\psi_2(f(X),f(Y))$.
\vskip 0.1truecm
\par
\noindent
Case 2: $\psi_2(X,Y)<0$.
\vskip 0.1truecm
\par
\noindent
Without loss of generality we may assume that $X=(0,0)$.
Let $Y=(a,b)$, $A=(-2ib, 2ia)$, $B=(2ib, -2ia)$.
Then $a^2+b^2=\psi_2(X,Y) \in \Q \cap (-\infty,0)$. Thus
$\psi_2(A,X)=\psi_2(B,X)=-4(a^2+b^2) \in \Q \cap (0,\infty)$,
$\psi_2(A,B)=-16(a^2+b^2) \in \Q \cap (0,\infty)$,
$\psi_2(A,Y)=\psi_2(B,Y)=-3(a^2+b^2) \in \Q \cap (0,\infty)$.
We show that
$$S_{XY}:=S_{AX} \cup S_{BX} \cup S_{AB} \cup S_{AY} \cup S_{BY}$$
satisfies the condition of the theorem;
the finite sets $S_{AX}$, $S_{BX}$, $S_{AB}$, $S_{AY}$, $S_{BY}$
are known to exist by the proof in case 1.
Assume that $f:S_{XY} \to {\K}^2$ preserves unit distance.
Then $f$ preserves the distances between $A$ and $X$, $B$ and $X$, $A$ and $B$, $A$ and $Y$, $B$ and $Y$.
By Proposition 3 the Cayley-Menger determinant
$\Delta(f(X),f(Y),f(A),f(B))$ equals $0$ i.e.
$$
\det \left[
\begin{array}{ccccc}
0 &      1      &      1      &       1      &       1      \\
1 &      0      &      t      &  -4(a^2+b^2) &  -4(a^2+b^2) \\
1 &      t      &      0      &  -3(a^2+b^2) &  -3(a^2+b^2) \\
1 & -4(a^2+b^2) & -3(a^2+b^2) &       0      & -16(a^2+b^2) \\
1 & -4(a^2+b^2) & -3(a^2+b^2) & -16(a^2+b^2) &       0      \\
\end{array}
\right]
=0
$$
where $t=\psi_2(f(X),f(Y))$. Computing this
determinant we obtain
$$
32(a^2+b^2)(t-(a^2+b^2))^2=0.
$$
Therefore $t=a^2+b^2=\psi_2(X,Y)$.
\vskip 0.1truecm
\par
\noindent
Case 3: $\psi_2(X,Y)=0$.
\vskip 0.1truecm
\par
\noindent
There exists a complex isometry of ${\C}^2$
that sends $(0,0)$ to $X$ and $(1,i)$ to $Y$, so without loss of
generality we may assume that $X=(0,0)$ and $Y=(1,i)$.
Let $A=(-1,0)$, $B=(1,0)$. Then
$\psi_2(A,X)=\psi_2(B,X)=1$, $\psi_2(A,B)=4$,
$\psi_2(A,Y)=3$, $\psi_2(B,Y)=-1$. We show that
$$
S_{XY}:=S_{AX} \cup S_{BX} \cup S_{AB} \cup S_{AY} \cup S_{BY}
$$
satisfies the condition of the theorem;
the finite sets $S_{AX}$, $S_{BX}$, $S_{AB}$, $S_{AY}$, $S_{BY}$
are known to exist by the proofs in cases 1 and 2.
Assume that $f:S_{XY} \to {\K}^2$ preserves unit distance.
Then $f$ preserves the distances between $A$ and $X$, $B$ and $X$, $A$ and $B$, $A$ and $Y$, $B$ and $Y$.
Since $1=\psi_2(f(A),f(X)) \neq \psi_2(f(A),f(Y))=3$
we conclude that $f(X) \neq f(Y)$.
By Proposition 3 the Cayley-Menger determinant
$\Delta(f(X),f(Y),f(A),f(B))$ equals $0$ i.e.
$$
\det \left[
\begin{array}{rrrrr}
~0 & ~1 &  1 & ~1 &  1 \\
~1 & ~0 &  t & ~1 &  1 \\
~1 & ~t &  0 & ~3 & -1 \\
~1 & ~1 &  3 & ~0 &  4 \\
~1 & ~1 & -1 & ~4 &  0 \\
\end{array}
\right]
=0
$$
where $t=\psi_2(f(X),f(Y))$. Computing this
determinant we obtain $-8t^2=0$.
Therefore $t=0=\psi_2(X,Y)$.
\vskip 0.1truecm
\par
\noindent
The proof is finished.
\vskip 0.1truecm
\par
As a corollary of Theorem 3 we get:
\vskip 0.1truecm
\par
\noindent
{\bf Theorem 4.} If $X,Y \in {\C}^2$, $\psi_2(X,Y) \in \Q$ and
$X \neq Y$, then any map $f:{\C}^2 \to {\K}^2$ that preserves unit distance
satisfies $\psi_2(X,Y)=\psi_2(f(X),f(Y))$ and $f(X) \neq f(Y)$.
\vskip 0.1truecm
\par
If $X,Y \in {\C}^2$ and $\psi_2(X,Y) \not\in \Q$, then there exists
$f:{\C}^2 \to {\C}^2$ that does not preserve the distance between
$X$ and $Y$ although $f$ satisfies
$\psi_2(S,T)=\psi_2(f(S),f(T))$ for all $S,T \in {\C}^2$
with $\psi_2(S,T) \in \Q$;
the proof is analogous to the proof of Theorem~1.
\vskip 0.1truecm
\par
\noindent
Using {\boldmath $\left(\bullet\right)$} the author proved (\cite{Tyszka2004b}):
\begin{description}
\item{{\boldmath $\left(\blacksquare\right)$}}
each unit-distance
preserving mapping $f:{\R}^2 \to {\F}^2$ has
a form $I \circ (\rho,\rho)$, where $\rho: \R \to \F$ is
a field homomorphism and $I: {\F}^2 \to {\F}^2$ is
an affine mapping with orthogonal linear part.
\end{description}
\vskip 0.1truecm
\par
\noindent
Using {\boldmath $\left(\blacksquare\right)$}
and Theorem 4 the author proved (\cite{Tyszka2005}):
\par
\noindent
each unit-distance preserving mapping $f:{\C}^2 \to {\C}^2$ has
a form $I \circ (\gamma,\gamma)$, where $\gamma: \C \to \C$ is
a field endomorphism and $I: {\C}^2 \to {\C}^2$ is
an affine mapping with orthogonal linear part.

Apoloniusz Tyszka\\
Technical Faculty\\
Hugo Ko\l{}\l{}\k{a}taj University\\
Balicka 116B, 30-149 Krak\'ow, Poland\\
E-mail address: {\it rttyszka@cyf-kr.edu.pl}
\end{document}